\documentclass{aptpub}
\usepackage{epsfig}

\authornames{Johnson Atkinson}
\shorttitle{Gambler's Ruin:  Some Aspects of Coin Flipping}

\usepackage{graphicx}
\usepackage{epstopdf}
\usepackage{epstopdf}

\catcode`\^^?=9 

\def\be{\begin{equation}}
\def\ee{\end{equation}}
\def\ba{\begin{eqnarray}}
\def\ea{\end{eqnarray}}
\def\nn{\nonumber}
\def\ban{\begin{eqnarray*}}
\def\ean{\end{eqnarray*}}
\def\mref#1{Eq.(\ref{eq:#1})}
\def\nref#1{(\ref{eq:#1})}
\def\mlab#1{\label{eq:#1}}

\begin{document}

\title{ Gambler's Ruin? \\  \phantom{i}{\small Some aspects of coin tossing}\\}

\authorone[Illinois Institute of Technology]{Porter W. Johnson}

\addressone{Illinois Institute of Technology, Chicago, IL 60616, U.S.A.}

\authortwo[University of Groningen]{David Atkinson}

\addresstwo{ University of Groningen, 9712 GL   Groningen, The Netherlands}

\begin{abstract}
What is the average number of coin tosses needed before a particular sequence of heads and tails first turns
up? This problem is solved in our paper, starting with doubles, a tail, followed by a head,
turns up on average after only four tosses, while six tosses are needed for two successive heads.
The method is extended to encompass the triples  head-tail-tail and head-head-tail, but head-tail-head
and head-head-head are surprisingly more recalcitrant. However, the general case is finally solved by
using a new algorithm, even for relatively long strings. It is shown that the average number of tosses
is always an even integer.
\end{abstract}



\vskip13mm

\section{Do heads like tails?}
The classic example of the gambler's fallacy in coin tossing is that a long run of heads
is more likely to be succeeded by a tail than by yet another head. 
A variant is as follows. A gambler must make a stake of 5 euros.
He can decide whether to bet on a head-head, or a tail-head outcome. In the first
case the gambler wins $n$ euros if $n$ tosses are made until two successive
heads first appear, and in the second case he wins $n$ euros if $n$ tosses
are made until the sequence tail-head first appears. Does
it make any difference which of the two wagers is favoured? 

The gambler might reason that, assuming the coin to be fair, an
equal number of heads and tails should come up in the long run, and
therefore it is more likely that the sequence tail-head will come up
than the  sequence head-head. So he should bet on head-head, since
it is likelier that more tosses will be needed than would be required
for tail-head. Is the gambler wrong again? Are tail-head and head-head
equally likely? 

Notwithstanding any first intuition one might have, it turns out that
the gambler is well-advised to bet on head-head, and not on tail-head.
Indeed, on average he would win one euro per game with head-head,
but lose one euro per game with tail-head. [1], [2]

In Section 2 we treat the problem of `doubles', calculating the average
number of tosses needed before a tail, followed by a head (TH) first comes up,
and then  the average number before two heads (HH) first come up. These
averages are computed by first obtaining the probabilities that TH or HH
first come up after $m$ tosses, and then calculating the mean value of $m$.
Since there is symmetry under an interchange of heads and tails,  TH and
HT yield the same average number of tosses, and similarly for TT  and HH.  
In the HH case the probability turns out to be related to a Fibonacci
number.

In Section 3 the probabilities associated with the doubles are recalculated
by using recurrence relations; and this method is then applied
to the  `triples' HTT and HHT. The corresponding calculations for HTH and HHH
prove cumbersome using the technique of this section, and they are not therefore
attempted.

The calculations in Sections 2 and 3 unexpectedly result in averages that are integers:
why should the average number of tosses be a whole number? The answer to this
question is given in Section 4, in which a  new method
is introduced that enables all the findings of the previous sections to be
recovered, as well as providing the averages pertaining to the remaining triples. This
technique makes it clear that all the averages must in fact be even integers.
At the end of the paper a table is presented of the average numbers of tosses,
not only for the doubles and triples, but also for the quadruples, quintuples
and sextuples, all of which can be readily treated by the new method.


\section{Doubles}
We wish  to calculate the probability that a tail
first comes up in the $(m-1)$st toss, followed
immediately by a head in the $m$th toss. Representing a tail by 0
and a head by 1,  and the sequence of tosses by a string of zeroes
and ones, we see that this is the same as the probability that
for the first time the
$(m-1)$st and $m$th positions are occupied by the sequence 01.
Suppose that there are $n$ tosses in all, so that we have
$2^n$ distinct configurations of zeroes and ones, each
being equiprobable. To calculate the probability of interest, we
first consider the number of all the distinct configurations in
which the sequence 01 is {\sf  absent}. If $n=1$, so that there is
only one toss, the sequence 01 is absent.
For two tosses,  $n=2$, three of the four configurations are
not 01. For $n=3$, the number of configurations lacking 01
is four, and for $n=4$ it is five.  Table 1  enables one to
see at a glance why, for general $n$, the number of configurations
that do not contain 01 is $n+1$. For to progress from $n$ to
$n+1$ one can add a 0 to the right of each $n$-sequence without
creating 01, and one can also add a 1 to the right of
the configuration containing only ones. No other configuration
tolerates the adjunction of 1, since  all the rest have a 0 in
the rightmost position. 
 
\begin{picture}(0,140)(-98,-50)
\put(20,55){0}
\put(20,40){1}
\put(50,55){00}
\put(50,40){10}
\put(50,25){11}
\put(80,55){000}
\put(80,40){100}
\put(80,25){110}
\put(80,10){111}
\put(110,55){0000}
\put(110,40){1000}
\put(110,25){1100}
\put(110,10){1110}
\put(110,-5){1111}
\put(-20,-30){Table 1: Sequences {\sf  not} containing 01}
\end{picture}

How many different sequences of length $n$  are there, such that
the $(m-1)$st and $m$th positions are occupied by  01, and such that
the first $m-2$ positions do not contain the sequence 01? We have
just deduced that there will be $(m-2)+1=m-1$ such sequences in a
string of length $m-2$. There are $n-m$ positions after the 01, and
it is not important if 01 recurs, so there are $2^{n-m} $ distinct
sequences here. Therefore the total number we are seeking is 
$(m-1)\, 2^{n-m} $. Since the total number of distinct sequences
is $2^{n} $, we must normalize by this to obtain the probability,
$p_m$, that  01 first occurs in the $(m-1)$st and $m$th positions
respectively:
\[
p_m(01)=(m-1)\, 2^{-m}
\] 
It is now straightforward to compute the average value of $m$, namely 
\[
N_{01}=\sum_{m=2}^\infty m\, p_m(01)=
\sum_{m=2}^\infty m\, (m-1)\, 2^{-m}=4\,.
\]
In the long run a tail, followed by a head, will come up after 4 tosses,
so if our gambler were to put down his stake of 5 euros on the tail-head
option, he would lose an average of one euro per game.

What difference does it make if we now require that a head comes up in
the $(m-1)$st toss, followed immediately by another head in the $m$th
toss? What is the probability that the $(m-1)$st and $m$th positions in
a string of zeroes and ones are occupied for the first time by the sequence 11?
We consider now the number of all the distinct configurations
in which the sequence 11 is absent. If $n = 1$, with possibilities, 0 or
1, the sequence 11 is absent. For two tosses, $n = 2$, three of the four configurations
are not 11. For $n = 3$, the number of configurations lacking 11 is five,
and for $n = 4$ it is eight. These configurations are shown in Table 2,
and evidently this array is quite different from that in Table 1.  

\begin{picture}(0,180)(-80,-90)
\put(20,55){0}
\put(20,40){1}
\put(50,55){00}
\put(50,40){10}
\put(50,25){01}
\put(80,55){000}
\put(80,40){100}
\put(80,25){010}
\put(80,10){001}
\put(80,-5){101}
\put(110,55){0000}
\put(110,40){1000}
\put(110,25){0100}
\put(110,10){0010}
\put(110,-5){1010}
\put(110,-20){0001}
\put(110,-35){1001}
\put(110,-50){0101}
\put(-20,-75){Table 2: Sequences {\sf  not} containing 11}
\end{picture}

This pattern yields a Fibonacci sequence, which is defined by the recursion
\[
F_0=0\hspace{16mm}F_1=1\hspace{16mm}F_n=F_{n-1}+F_{n-2}\,,
\] 
so $F_2=1$, $F_3=2$, $F_4=3$, $F_5=5$, $F_6=8$, and so on. 
The number of configurations  corresponding to $n$ bits is $F_{n+2}$.

We should ask how many different sequences of length $n$ there are,
such that  the $(m-1)$st and $m$th positions are occupied for the first
time  by  11. We have just shown that there are $F_m$ different
sequences of length $m-2$ that lack 11. But this is not sufficient,
for some of these $F_m$ sequences end with 1, which is disallowed,
since we want 11 to be in  the $(m-1)$st and $m$th positions, not the
$(m-2)$nd and $(m-1)$st positions. So we must count only those of the
$F_m$ sequences that end in 0. There are precisely $F_{m-1}$ of
these sequences, as can be seen from Table 2. Indeed, they are
the sequences of $m-3$ bits that lack 11, augmented by a 0 at
the right end. 

The total number of states in a sequence of $n$ bits that have
11 for the first time in the $(m-1)$st and $m$th positions is
therefore  $F_{m-1}\, 2^{n-m} $. On normalizing this we obtain
the probability, $p_m(11)$, that  11 first occurs in the $(m-1)$st
and $m$th positions respectively:
\[
p_m(11)=F_{m-1}\, 2^{-m}
\] 
The average value of $m$ is therefore
\be
N_{11}=\sum_{m=2}^\infty m\, p_m(11)=
\sum_{m=2}^\infty m\, F_{m-1}\, 2^{-m}=6\,.\mlab{N11}
\ee
The above sum was evaluated from the definition of Fibonacci numbers.
In the long run a head, followed by another head, will come up after
6 tosses, so if our gambler were to put down his stake of 5 euros
on the head-head option, he would gain an average of one euro per game.

\section{Recursive solution}
We now repeat the above calculations of $N_{01}$ and $N_{11}$ with a notation
that leads directly to recursion relations. The new method will then be generalized to
sequences involving 3  consecutive bits (triples). Let $K_m(01)$ be the number of different strings
of zeroes and ones of length $m$, such  that  01 occurs for the only time at the right-hand end of
the string, as in the previous section. Similarly,  $K_m(11)$ is the number of different strings
of length $m$, such  that  11 occurs for the only time at the right-hand end of that string.

Consider first 01, and define $M_m(01;j)$ to be the number of  different
strings of length $m$, such  that  01 does {\sf  not} occur, and such that the rightmost bit is
$j$ ($j$ = 0 or 1). We expect a recursion relation of the sort 
\[
M_m(01;j)=a\, M_{m-1}(01;0)+b\, M_{m-1}(01;1)\,,
\]
for $M_{m-1}(01;j)$ is, like $M_m(01;j)$, free of the sequence 01, and it is one bit shorter.
Here $a$ and $b$ are 0 or 1, and we find 
\ba
M_m(01;0)&=&\, M_{m-1}(01;0)+\, M_{m-1}(01;1)
\nn  \\  
M_m(01;1)&=&\,  M_{m-1}(01;1)\mlab{iter1}
\ea
Note that $M_{m-1}(01;0)$ is missing in the second equation, i.e. $a=0$, for we may not
add a 1 at the end of a string of length $m-1$ that terminates with a 0, because that would
produce a string which terminates with 01, contrary to the definition of $M_m(01;1)$.
The second line of \nref{iter1} says that $M_{m}(01;1)$ is independent of $m$, and since
$M_1(01;1)$ is clearly one, we find  $M_{m}(01;1)=1$, and the first  line of \nref{iter1}  becomes 
\ban
M_m(01;0)&=&\, M_{m-1}(01;0)+1\,.
\ean
The solution is $M_m(01;0)=m$,  since $M_1(01;0)$ is also clearly one.  A string of length $m$
that contains 01 only at its right-hand end can be thought of as a string of length $m-1$ that
does not contain 01, but which has 0 at {\sf  its}  right-hand end, plus a 1 at the extreme right. This means that 
\be
K_m(01)=M_{m-1}(01;0)=m-1\,,
\ee
which agrees with our earlier result.

Let us now consider 11 and similarly define $M_m(11;j)$ to be the number of  different
strings of length $m$, such  that  11 does  not  occur, and such that the rightmost bit
is 0 or 1. The analogue of \mref{iter1} is
 \ba
M_m(11;0)&=&\, M_{m-1}(11;0)+\, M_{m-1}(11;1)
\nn  \\  
M_m(11;1)&=&\,  M_{m-1}(11;0)\mlab{iter2}
\ea
Here it is $M_{m-1}(11;1)$ that is missing, for adding a 1 to the end of that sequence would produce
a string containing 11. Replace $m$ by $m-1$ in the second line of \nref{iter2} and
substitute for $M_{m-1}(11;1)$ in the first line of  \nref{iter2}:
 \ban
M_m(11;0)&=&\, M_{m-1}(11;0)+\, M_{m-2}(11;0)\,.
\ean
This is the Fibonacci recursion, and since $M_1(11;0)=1$ and $M_2(11;0)=2$, it follows that 
$M_m(11;0)=F_{m+1}$. Clearly,
\[
K_m(11)=M_{m-1}(11;1)=M_{m-2}(11;0)=F_{m-1}\,,
\]
once more agreeing with the result of Section 2. 

We now break new ground by considering triples. How many strings of length $m$ are there for
which 100 occurs only at the right-hand end? We call this number $K_m(100)$. Let  $M_m(100; ij)$
be the number of strings of length $m$ that do not contain the sequence 100, and which end
with $ij$, i.e. either 00, 01, 10 or 11. The relevant recursion relations are now 
 \ba
M_m(100;00)&=& M_{m-1}(100;00)
\nn  \\  
M_m(100;01)&=& M_{m-1}(100;00)+M_{m-1}(100;10)
\nn  \\  
M_m(100;10)&=& M_{m-1}(100;01)+M_{m-1}(100;11)
\nn  \\  
M_m(100;11)&=&\,  M_{m-1}(100;01)+M_{m-1}(100;11) \mlab{iter3}
\ea
The last three lines have two terms on the right, corresponding to the possibility of adding
the relevant bit to the right of the string of length $m-1$. The first line lacks $M_{m-1}(100;10)$,
because the adjunction of a 0 to the right would produce a string ending in 100. This line shows
that $M_m(100;00)$ is independent of $m$, so $M_m(100;00)=M_2(100;00)=1$. 
The second line of  \mref{iter3} is therefore 
\be
M_m(100;01)=1+M_{m-1}(100;10)\mlab{iter3a}\,.
\ee
We next note that the right-hand sides of the third and fourth lines are identical,
so $M_m(100;10)=M_m(100;11)$, and therefore 
\be
M_m(100;10)=M_{m-1}(100;01)+M_{m-1}(100;10)  \mlab{iter3b}\,.
\ee
Change $m$ into $m-1$ in \nref{iter3a}, and substitute the result into  \nref{iter3b}:
\be
M_m(100;10)=1+M_{m-2}(100;10)+M_{m-1}(100;10)  \mlab{iter3c}\,.
\ee
It is clear that $1+M_m(100;10)$ satisfies the Fibonacci recurrence relation. Since  
 $M_3(100;10)=2$, we conclude that 
$M_m(100;10)=F_{m+1}-1$. Thus the number of strings ending in 100 is 
\be
K_m(100)=M_{m-1}(100;10)=F_m-1\,,\mlab{K100}
\ee
so the average value of $m$ is 
\be
N_{100}=
\sum_{m=3}^\infty m\, [F_m-1]\, 2^{-m}=8\,.\mlab{N100}
\ee

We next  consider the triple 110.  Let  $M_m(110; ij)$ be the number of strings of length $m$ that do
not contain the sequence 110, and which end with $ij$, i.e. either 00, 01, 10 or 11. The relevant recursion
relations are now 
 \ba
M_m(110;00)&=& M_{m-1}(110;00)+M_{m-1}(110;10)
\nn  \\  
M_m(110;01)&=& M_{m-1}(110;00)+M_{m-1}(110;10)
\nn  \\  
M_m(110;10)&=& M_{m-1}(110;01)
\nn  \\  
M_m(110;11)&=&\,  M_{m-1}(110;01)+M_{m-1}(110;11) \mlab{iter4}
\ea
The third line lacks $M_{m-1}(110;11)$, because the adjunction of a 0 to the right would produce a string containing 110. 
The right-hand sides of the first and second lines are identical, so $M_m(110;00)=M_m(110;01)$, and therefore 
\be
M_m(110;01)=M_{m-1}(110;01)+M_{m-2}(110;01)\mlab{iter4h}
\ee
where the third line has also been used. This is the Fibonacci recursion again, and since 
$M_3(110;01)=2$ it follows that $M_m(110;01)=M_m(110;00)=F_m$, and also that $M_m(110;10)=F_{m-1}$.
The fourth line in \mref{iter4} becomes
$
M_m(110;11)=  F_{m-1}+M_{m-1}(110;11)\,,
$
and by iteration this yields 
\be
M_m(110;11)= \sum_{p=1}^{m-1}F_p=F_{m+1}-1
\ee
Finally, the number of strings ending in 110 is 
\[
K_m(110)=M_{m-1}(100;11)=F_m-1\,,
\]
From \mref{K100} we see that this is the same formula as the one we derived for the triple 100, so we also find  
$N_{110}=8$.

The above method turns out to be inadequate for the calculation of $N_{101}$ and $N_{111}$:
an analogous treatment would involve more complicated numbers than those of  Fibonacci, and so
in the following section we examine a novel approach that turns out to be rather simpler.

\section{All averages are even integers}
An interesting feature of the average numbers, 4, 6 and 8 that we have calculated so far is that they are integral. 
The number of tosses in any one game is of course an integer; but why should their {\sf  average} be a whole number too? 
In this section we shall demonstrate that the average
number of tosses needed before a given string, 
\be
T_m=(t_1t_2\ldots t_m)\mlab{tm1}\,,
\ee
comes up for the first time is in general an integer, and moreover an {\sf  even}
integer.  Consider first a string of length $n$, 
\be
S_n=(s_1s_2\ldots s_n)\mlab{tm3}\,,
\ee
that does {\sf  not} contain $T_m$ as a substring; and let $\sigma_n$ be the number of such strings. 
A string of length $m+n$ that contains  $T_m$ just once, right at the end, has the form 
\be
(S_nT_m)=(s_1s_2\ldots s_n t_1t_2\ldots t_m)\,. \mlab{SnTmgen}
\ee
Let $\tau_{m+n}$ be the number of such strings of length $m+n$ that terminate with $T_m$,
but which do not contain any other instance of $T_m$.\footnote{In the notation of
Section 3, $\tau_{m+n}$ corresponds to $K_{m+n}(T_m)$,  and  $\sigma_n$  
corresponds to $M_n(T_m;0)+M_n(T_m;1)$. }
 
To make a string of the type $S_n$ from $S_{n-1}$ one can add either $s_n=0$ or $s_n=1$
at the end, and thus there will often be twice as many strings of the type
$S_n$ as of the type $S_{n-1}$. However, if  the addition of 0 or 1 to a string   $S_{n-1}$
results in the completion of the string $T_m$, the new  string should not be counted as one
of the tally $\sigma_n$, but rather as one of the tally $\tau_n$, since it has become a
string of length $n$, of which the last $m$ bits constitute $T_m$. Thus 

 \be
2\sigma_{n-1}= \sigma_n +\tau_n \, \mlab{sigiter} \,.
 \ee
If $m=1$, $ \sigma_1=1$ and $\tau_1=1$, but if $m\ge 2$, $ \sigma_1=2$ and $\tau_1=0$,
so in all cases \nref{sigiter} is  true for $n=1$ if we formally define  $ \sigma_0=1$.
With this understanding,  \nref{sigiter} is valid for $n=1,2,3\ldots$
 
By iteration we find from \mref{sigiter} that 
 \[
  \sigma_n=-\tau_n-2\tau_{n-1}-\ldots -2^{n-m}\tau_m+2^{n-m+1}\sigma_{m-1}\,.
 \] 
Since one cannot obtain $T_m$ with only $m-1$ bits, 
$\sigma_{m-1}$ is equal to the total number of strings of length $m-1$, namely $2^{m-1}$, and thus 
 \be
  \sigma_n=2^n\left( 1-\sum_{j=m}^n \frac{\tau_j}{2^j}\right)
  =2^n\sum_{j=n+1}^\infty \frac{\tau_j}{2^j}\,,  \mlab{summation}
 \ee
where in the last step we have used the identity 
$\sum_{j=m}^\infty {\tau_j}/{2^j}=1$, which is proved in the Appendix --- \mref{sigiter3}. 

Dividing both sides of \mref{summation} by $2^n$ and summing, we find
\be
\sum_{n=0}^\infty \frac{\sigma_n}{2^n}= \sum_{n=0}^\infty\sum_{j=n+1}^\infty \frac{\tau_j}{2^j}\,.
\mlab{sigsum}
\ee
On changing the order of the summations, we can rewrite the  right-hand side in the form
 \be
 \sum_{j=0}^\infty \frac{\tau_j}{2^j}\sum_{n=0}^{j-1} 1
=\sum_{j=m}^\infty j \, \frac{\tau_j}{2^j}
\,.\mlab{sigsumsum}
 \ee
This is the average value of the lengths of the strings that terminate in $T_m$,
for $\tau_j$ is the number of them with length $j$, while $2^{-j}$ is their probabilistic weight. 
The lower limit of the summation on the right-hand side has been changed to $j=m$, since $\tau_j=0$ for $j<m$. 
In the notation of the previous sections, the sum \nref{sigsumsum} is $N_{T_m}$, so \mref{sigsum} becomes 
 \be
N_{T_m}=\sum_{n=0}^\infty \frac{\sigma_n}{2^n}\,, \mlab{Ntmend}
 \ee
where we recall that $\sigma_0=1$.

Let us first consider  the case in which $t_1=1$ and $t_j=0$ for $2\le j\le m$,
so the string of length $m$ has the form 
\[
T_m=(100\ldots 0)\mlab{tm22}\,.
\]
A string of length $m+n$ that contains  $T_m$ just once, right at the end, has the form 
\be
(S_nT_m)=(s_1s_2\ldots s_n100\ldots 0)\,. \mlab{SnTm}
\ee
The number of strings \nref{SnTm} is equal to the number of strings \nref{tm3}, there
being only one way to make a string $S_nT_m$  from a string  $S_n$, namely by
adding $T_m$ to its end, so 
 \be
 \sigma_n=\tau_{m+n}\,. \mlab{end}
 \ee
On substituting  \mref{end} into \mref{Ntmend}, we obtain
\[
N_{100\ldots 0}=\sum_{n=0}^\infty \frac{\tau_{m+n}}{2^n}= 2^m\sum_{j=m}^\infty \frac{\tau_j}{2^j}=2^m\,,
\]
where \mref{sigiter3} has  been used. 

The next case we shall treat is that in which $t_1=1$ for $1\le j\le m$, so the string of
length $m$ contains only ones:
\[
T_m=(111\ldots 1)\mlab{tm23}
\]
Consider the string of length $m+n$ of the form 
\be
(S_nT_m)=(s_1s_2\ldots s_n111\ldots 1)\,. \mlab{SnTm2}
\ee
The number of such strings is equal to the number of strings $(S_n)$, but this 
is not now equal to $\tau_{m+n}$, since some of the strings $S_n$ will end in 01, and
when $T_m$ is appended to those strings, they will belong to the tally $\tau_{m+n-1}$. Similarly,
some of the strings $S_n$ will end in 011, and when $T_m$ is appended to those strings,
they will belong to the tally $\tau_{m+n-2}$, and so on, up to the strings $S_n$ that end in zero,
followed by $m-1$ ones. When $T_m$ is appended to these strings, they will belong to the tally $\tau_{n+1}$.
So instead of \mref{end}, we have in the present case 
\be
\sigma_n=\sum_{j=1}^m \tau_{j+n}\,. \mlab{end2}
 \ee
Substituting  \mref{end2} into \mref{Ntmend}, we now find 
\[
N_{111\ldots 1}=\sum_{n=0}^\infty 2^{-n}\sum_{j=1}^m \tau_{j+n}= \sum_{j=1}^m 2^j\sum_{i=m}^\infty \frac{\tau_i}{2^i}=
\sum_{j=1}^m 2^j\
\,,
\]
where \mref{sigiter3} has again been used.  

$N_{100\ldots 0}$ and $N_{111\ldots 1}$ are two extreme cases. In general a string $T_m$ of length $m$ 
will generate a series 
\be
\sigma_n=\sum_{j=1}^m c_j\, \tau_{j+n}\,,\mlab{end3}
 \ee
 where always $c_m=1$, but where some of the other coefficients $c_j$ may be zero,
the others being equal to one. In this general case we find
 \be
 N_{T_m}=\sum_{n=0}^\infty 2^{-n}\sum_{j=1}^m c_j\,\tau_{j+n}=
\sum_{j=1}^m c_j\,2^j\sum_{i=m}^\infty \frac{\tau_i}{2^i}=
\sum_{j=1}^m c_j\,2^j \mlab{Ntmform}
\,,
 \ee
once more with use of \mref{sigiter3}. This expression shows
that $ N_{T_m}$ is always an even integer, being the sum of positive integral powers of two.
The coefficient $c_j$ is equal to zero if the string $(S_nT_m)$, truncated at order $n+j$, does
not  contain $T_m$ for its last $m$ bits; otherwise it is equal to one.  For a given $m$, the
following inequalities hold:
\[
 2^m\le  N_{T_m}\le 
\sum_{j=1}^m \,2^j =2^{m+1}-2\,,
\]
with  \mref{end} and \mref{end2}, respectively realizing the extreme possibilities.

By way of example, we shall show how to extract the coefficients $c_j$ of \mref{end3}
in a particular case in which $m=5$, namely the string 
$T_{10101}\,.$
On placing $S_n$ before $T_{10101}$ we obtain a string 
\be
(S_nT_5)=(s_1s_2\ldots s_n10101)\,,\mlab{SnT5}
\ee
and the number of such strings is the same as the number of strings of the
type $S_n$, namely $\sigma_n$. However, the 
strings \nref{SnT5} encompass not only those of length $5+n$
that terminate with $T_{10101}$, but also strings of length $5+n-2= 3+n$, for which the
first $10$ of $10101$ belongs to $S_n$, i.e. strings of the sort 
\[
(s_1s_2\ldots s_{n-2}1010101)\,.
\]
Similarly, \nref{SnT5} also includes strings of length $5+n-4 =1+n$ for
which the first $1010$ of $10101$ belongs to $S_n$.  
No such string could have terminated at order
$4+n$, since the first four bits of $10101$, namely $1010$, are not consistent with
the last four bits of the required termination; namely $0101$.  That is, 
\[
(s_1s_2\ldots s_{n-1}s_n10101)\,,
\]
cannot belong to the tally $\tau_{n+4}$, irrespective of whether $s_n$ is 0 or 1, so it
must be part of the tally $\tau_{n+5}$.
For a similar reason, no such string could have terminated at order $2+n$,
since the $(2+n)$th bit would be 0, whereas the last bit of the termination should be 1.
So we obtain
\[
\sigma_n=\tau_{n+1}+\tau_{n+3}+\tau_{n+5}\,,
\]
and thus
\[
N_{10101}=2+2^3+2^5=42
\,.
\]
The results shown in Table 3 were obtained by an application of this method.

\section*{Appendix: Convergence and completeness}
Divide both sides of \mref{sigiter} by $2^n$ and sum:
 \[
 \sum_{n=m}^p\frac{\sigma_{n-1}}{2^{n-1}}=\sum_{n=m}^p\frac{\tau_n}{2^n} 
 +\sum_{n=m}^p\frac{\sigma_n}{2^n} \mlab{sigiter2}\,,
 \]
  and so 
  \be
\sum_{n=m}^p\frac{\tau_n}{2^n} =
 \sum_{n=m-1}^{p-1}\frac{\sigma_n}{2^n}-\sum_{n=m}^p\frac{\sigma_n}{2^n}=
\frac{\sigma_{m-1}}{2^{m-1}}-\frac{\sigma_{p}}{2^{p}}
 \mlab{newsum}\,.
 \ee
 Now $\sigma_{m-1}=2^{m-1}$, since all $2^{m-1}$ strings of length $m-1$ are of the type $S_{m-1}$. It follows that 
\[
 \sum_{n=m}^p\frac{\tau_n}{2^n}\le 1\,,
\] 
 for all $p$, which means that the infinite series 
 $
  \sum_{n=m}^\infty{\tau_n}{/2^n}
 $
 is convergent, since all the terms are positive and the partial sums are uniformly bounded. 
 This convergence means in particular that ${\tau_n}{/2^n}$ tends to 0 in the limit that $n$ goes to infinity. 
However, from \mref{end3}  we see that 
 \[
\frac{\sigma_n}{2^n} =\sum_{j=1}^m c_j\, 2^j\frac{\tau_{j+n}}{2^{j+n}}
\]
goes to zero as $n$ tends to infinity at constant $m$.  
 Therefore from \mref{newsum} it follows that 
 \be
 \sum_{n=m}^\infty\frac{\tau_n}{2^n}=1-\lim_{p\rightarrow\infty}\frac{\sigma_{p}}{2^{p}}=1
 \mlab{sigiter3}\,.
 \ee
\vspace{7mm}

\noindent
{\bf Acknowledgement:} We would like to thank Barteld Kooi for having brought this interesting problem to our notice.

\newpage

 \[
\begin{array}{||c|c|l||}
\hline
\mbox{Num bits} &\mbox{Average} &  \mbox{Terminating String(s)} \\
\hline
& 4 &  {10} \\
2 & 6  & {11} \\   \hline

& 8  & {100}, \, {110} \\
3 & 10  & {101} \\
& 14  & {111}  \\  \hline

& 16  & {1000}, \, {1100}, \, {1110} \\
& 18  & {1001},  \, {1011}, \, {1101} \\
4 & 20  & {1010} \\
& 30  & {1111} \\   \hline
 
& 32 &   {10000}, \, {10100}, \, {11000}, \, {11010}, \, {11100}, \, {11110} \\
& 34 &   {10001}, \, {10011}, \, {10111},  \, {11001},  \, {11101} \\
& 36 &   {10010}, \, {10110} \\
5 & 38 &   {11011} \\ 
& 42 &   {10101} \\ 
& 62 &   {11111} \\ \hline

& 64  & {100000}, \, {101000}, \, {101100}, \, {110000}, \, {110010}, \\
& & {110100}, \, {111000}, \, {111010}, \, {111100}, \, {111110} \\
\cline{2-3}
& 66  & {100001}, \, {100011}, \, {100101}, \, {100111}, \, {101001}, \, {101011}, \\
& & {101111}, \, {110001}, \, {110101}, \, {111001}, \, {111101} \\
\cline{2-3}
6 & 68  & {100010}, \, {100110}, \, {101110} \\
& 70  & {110011}, \, {110111}, \, {111011} \\
& 72  & {100100}, \, {110110}\\
& 74  & {101101} \\
& 84  & {101010} \\
& 126  & {111111} \\ \hline
\end{array}
\]

\begin{center}
{\small\sf Average number of tosses for doubles, triples, quadruples, quintuples and sextuples}
\vspace{5mm}

{\bf Table 3}
\end{center}

\end{document}